\documentclass{article}

\newtheorem{theorem}{Theorem}[section]
\newtheorem{lemma}[theorem]{Lemma}
\newtheorem{corollary}[theorem]{Corollary}

\newcommand\qed{\begin{flushright} {\bf q.e.d.} \end{flushright} }
\newcommand\prf{\noindent {\bf Proof :}  }

\newcommand\nn{{\{0,1\}^n}}
\newcommand\mm{{\{0,1\}^m}}
\newcommand\at{ \{0, 1\}^t }

\newcommand\dfc{{\mbox{Def}_C^{n,s}}}
\newcommand\dfz{{\mbox{Def}_C^{0,s}}}
\newcommand\gn{\Gamma(n,s,k)}
\newcommand\gz{\Gamma(0,s,k)}
\newcommand\ig{{i\Gamma}}
\newcommand\igm{{\ig(m)}}

\newcommand\xx{{\mathbf x}}
\newcommand\yy{{\mathbf y}}
\newcommand\pp{{\mathbf p}}
\newcommand\qq{{\mathbf q}}
\newcommand\zz{{\mathbf z}}
\newcommand\ab{{\mathbf a}}
\newcommand\bb{{\mathbf b}}
\newcommand\vv{{\mathbf v}}
\newcommand\rw{\mbox{R}^w}

\begin{document}

\title{Consistency of circuit evaluation, extended resolution and total NP search problems\\
({\sf preliminary version})}

\author{Jan Kraj\'{\i}\v{c}ek}

\date{Faculty of Mathematics and Physics\\
Charles University in Prague}

\maketitle

\begin{abstract}
We consider sets $\gn$ of narrow clauses expressing that no definition of a size $s$ circuit
with $n$ inputs is refutable in resolution R in $k$ steps. We show that every CNF shortly refutable
in Extended R, ER, can be easily reduced to an instance of $\gz$ (with $s,k$ depending on the size of the 
ER-refutation) and, in particular, that
$\gz$ when interpreted as a 
relativized NP search problem is complete among all such problems
provably total in bounded 
arithmetic theory $V^1_1$. 

We use the ideas of implicit proofs from \cite{Kra-implicit,Kra-hard}
to define from $\gz$ a non-relativized NP search problem $\ig$ and we show 
that it is complete among all such problems 
provably total in bounded arithmetic theory $V^1_2$. The reductions are definable in $S^1_2$.
 
We indicate how similar results can be proved for some other propositional
proof systems and bounded arithmetic theories
and how the construction can be used to define specific random unsatisfiable formulas,
and we formulate two open problems about them.
\end{abstract}

Let $C$ be a size $s$ circuit with $n$ Boolean 
inputs $\xx = x_1, \dots, x_n$ and in the basis $0,1, \neg, \vee, \wedge$.
It is defined by $s$ instructions how to compute Boolean
values $\yy = y_1, \dots, y_s$, all of which have one of the
following forms:
\begin{itemize}

\item $y_i := x_u$ for some $u \le n$,

\item $y_i := 0$ or $y_i := 1$, 

\item $y_i := \neg y_j$ for some $j < i$,

\item $y_i := y_j \vee y_k$ or $y_i := y_j \wedge y_k$ for some $j,k < i$. 

\end{itemize}
The value of $y_s$ is the output value of $C$ and is denoted also as $C(\xx)$.
Let $\dfc(\xx,\yy)$ be the canonical 3CNF formula expressing the conjunction of all instructions.
For example, instruction $y_i := 0$ is represented by one clause $\{y_i^0\}$,
$y_i := \neg y_j$ by 2 clauses $\{y_i, y_j\}, \{y_i^0, y_j^0\}$ and
instruction $y_i := y_j \vee y_k$ is represented by three clauses
$$
\{y_j^0, y_i^1\}\ ,\ \{y_k^0, y_i^1\}\ ,\ \{y_i^0, y_j^1, y_k^1\}\ 
$$
where for a literal $\ell$ define $\ell^1 := \ell$ and $\ell^0 := \neg \ell$.
$\dfc(\xx,\yy)$ has at most $3s$ clauses.

It is easy to prove in (propositional) resolution proof system R 
that the computation of $C$ is unique: in $O(i)$ steps derive from
$\dfc(\xx,\yy) \cup \dfc(\xx,\zz)$ clauses $\{y_i^1, z_i^0\}, \{y_i^0, z_i^1\}$ expressing
that $y_i \equiv z_i$. The whole proof of $y_s \equiv z_s$ has $O(s)$ clauses and
its structure is quite close to that of $C$.

But can we prove equally easily that a computation of $C$ on $\xx$ exists? This question is
in propositional logic represented by the question whether 
$\dfc(\xx,\yy)$ is consistent, i.e. not refutable, and we take as our refutation system
R (more precisely, its slight technical variant $\rw$ defined in Section \ref{1}). 
Given $n \geq 0, s, k \geq 1$ we shall define a set 
$\gn$ of narrow clauses such that satisfying assignments for $\gn$
would be precisely $k$ step $\rw$-refutations of sets $\dfc(x,y)$. 
Our question can be then phrased as follows:
How hard it is to refute $\gn$?

We will, in fact, concentrate on the case $n=0$ in which sets $\gz$ talk about refutations
of $\dfz(\yy)$, sets of clauses defining a straight-line program $C$
computing Boolean constants (i.e. $C$
has no inputs $\xx$).
Using standard techniques of proof complexity we 
show that sets $\gz$ express the reflection principle for Tseitin's \cite{Tse68}
Extended resolution ER, 
and hence any proof system that refutes these sets by polynomial size proofs
has at most polynomial slow-down over ER (it simulates it in the standard terminology). 
In fact, due to the combinatorial transparency of $\gz$
we use rather only the idea how reflection principles work rather than any
"technique" surrounding them. Further, 
the simulation yields straightforwardly
a reduction of unsatisfiable CNFs $\Delta$ 
to $\gz$ where $k$ depends on the size of an ER-refutation of $\Delta$, 
if we interpret them as relativized total NP search problems
with oracles representing truth assignments.

We will also show, using the idea of implicit proofs 
from \cite{Kra-implicit,Kra-hard}, how to define ordinary
(i.e. non-relativized) total NP search problems 
$\ig$ and show that these are complete among 
all NP search problems provably total
in theory $V^1_2$ of Buss \cite{Buss-book}. The reductions are definable
in $S^1_2$.
Another total NP search problems with this property have been defined earlier by
Kolodziejczyk, Nguyen and Thapen \cite{KNT} and recently by Beckmann and Buss \cite{Bec-Bus}.

We shall conclude the with remarks how to modify the construction for some other proof systems
and how to use it to define random unsatisfiable formulas, and we formulate two open problems.

Background from proof complexity we assume is standard and
can be found in \cite{kniha}. Only Section \ref{22.8.15a}
presupposes knowledge of a specific material from \cite{Kra-implicit,Kra-hard}; we explain there
the underlying ideas and give 
precise references but we shall not repeat that material here.

\section{Formalization: sets $\gn$} \label{1}

We shall first augment R a bit to make it technically more convenient. First, we shall
allow also constant $1$ in clauses and allow as new initial clauses all $C$ containing $1$
(we shall call these new initial clauses {\em 1-axioms}). Second, 
we add {\em the weakening rule}:
$$
\frac{C}{D}\ ,\ \mbox{ if }\ C \subseteq D\ .
$$ 
Denote this augmented resolution system $\rw$. The reason for the modifications
is that one can substitute constants for variables in an $\rw$-proof and it remains an
$\rw$-proof (delete all literals evaluated to $0$ and replace resolution inferences on
variables substituted for by weakenings). Additional reason for the weakening rule is that 
otherwise it is a bit cumbersome to talk about a derivation of $D$ from $C \subseteq D$:
as R is a refutation system one has to talk instead of refuting the set of clauses
$$
\{ C\} \cup \{ \{\ell^0\}\ |\ \ell \in D\}\ 
$$
and such derivations is even more cumbersome to concatenate (recall from
the introduction that $\ell^0 :=\neg \ell$).

\bigskip

Fix $n \geq 0$ and $s,k \geq 1$. Formula $\gn$ talks about a potential $k$-step 
$\rw$-refutation
of $\dfc$ for an unspecified $C$ (it is coded by atoms of $\gn$). 
For the purpose of the following discussion call these steps
$D_1, \dots, D_k$.

Clauses $D_i$ may contain constant $1$ or literals corresponding to $\xx, \yy$ variables,
i.e. all together up to $1 + 2(n+s)$ different objects. Formula $\gn$ will thus
use
\begin{itemize}

\item atoms $q^u_i$ with $u = 1, \dots, k$ and 
$i \in \{-(n+s), \dots, -1, 0, 1, \dots, (n+s)\}$

\end{itemize}
The intended
meaning of these is:
\begin{itemize}

\item $q^u_0 = 1$ iff $1 \in D_u$,

\item $q^u_i = 1$ for $i = 1, \dots, n$ iff $x_i \in D^u$, and 
for $i = -1, \dots, -n$ iff $x^0_i \in D^u$,

\item $q^u_{n+j} = 1$ for $j = 1, \dots, s$ iff $y_j \in D^u$, and 

$q^u_{-n+j} = 1$
for $j = -1, \dots, -s$ iff $y^0_j \in D^u$.

\end{itemize}
We shall call these variables $\qq$-variables and their set $\qq$.

There will be also $\pp$-variables $p_{u,v}$, $u = 1, \dots, k$ and $v = 1, \dots, t$
(we shall specify $t$ in a moment). The intended meaning is that an assignment 
$\ab^u \in \at$ for $\pp_u = p_{u,1}, \dots, p_{u,t}$ uniquely determines complete
information about how $D_u$ was inferred from earlier clauses and
if $D_u \in \dfc$ it also contains information assuring that $\dfc$ clauses have the
right form. To simplify the notation we shall assume that $k \geq 3s$ and $s > n$ and that the clauses
of $\dfc$ are listed as first $3s$ clauses $D_1, \dots, D_{3s}$, with 
$D_{3r-2}, D_{3r-1}, D_{3r}$ defining the instruction for $y_r$ (if the instruction
needs only one or two clauses the other are dummy, say $\{1\}$).

There are at most $2 + n + (r-1) + 2(r-1)^2 \le O(k^2)$ instructions how to compute $y_r$ 
and $\ab^u$ has to specify this uniquely for $u = 1, \dots, 3s$. For $u = 3s+1, \dots, k$
we need $\ab^u$ to specify by which rule and from which earlier clauses was $D_u$
inferred: there are at most $2 + (u-1) + (2+n+s)(u-1)^2 \le
k^3$ possibilities. Thus if we pick $t := 3 \log k$, $\at$ has enough room to 
encode by its elements all possible situations.  

It will be convenient to describe the clauses forming $\gn$ as sequents
$$
\ell_1, \dots, \ell_e\ \rightarrow\ \ell_{e+1}, \dots, \ell_f
$$
representing the clause
$$
\ell^0_1, \dots, \ell^0_e, \ell_{e+1}, \dots, \ell_f\ .
$$
For $\ab \in \at$ let $\pp_u(\ab)$ be the set of literals
$$
(p_{u,1})^{a_1}, \dots, (p_{u,t})^{a_t}\ .
$$
That is, $\ab$ is the unique truth assignment satisfying the conjunction of
literals in $\pp_u(\ab)$.

\bigskip

The set $\gn$ consists of the following clauses divided into five groups:

\begin{enumerate}

\item[$\gamma$1.] For $u \in \{3r-2, 3r-1, 3r\}$ for $r = 1, \dots, s$,
if $\ab \in \at$ does not specify a valid instruction for computing
$y_r$ then $\gn$ contains clause
$$
\pp_u(\ab) \rightarrow\ .
$$

\item[$\gamma$2.] For $u \in \{3r-2, 3r-1, 3r\}$ for $r = 1, \dots, s$,
if $\ab \in \at$ does specify a valid instruction for computing
$y_r$ then we know about constant $1$ and about every $\xx$- and $\yy$-variable
whether or not it occurs in $D_u$ and whether or not this occurrence is positive
or negative.  Hence we include in $\gn$ for every $\qq$-variable $q^u_i$
exactly one of the clauses
$$
\pp_u(\ab) \rightarrow q^u_i\ \mbox{ or }\ 
\pp_u(\ab) \rightarrow \neg q^u_i 
$$
as specified by $\ab$.

\item[$\gamma$3.] For $u = 3s+1, \dots, k$, if $\ab \in \at$ does not specify a valid inference
for $D_u$, $\gn$ contains clause
$$
\pp_u(\ab) \rightarrow\ .
$$

\item[$\gamma$4.] For $u = 3s+1, \dots, k$, if $\ab \in \at$ does specify a valid inference
for $D_u$, three cases can happen:

\begin{enumerate}

\item $D_u$ was inferred from $D_v, D_w$ resolving literal $\ell$, where $\ell \in D_v$
and $\ell^0 \in D_w$, and $\ell$ an $\xx$- or an $\yy$-literal. 

Let 
$i \in \{(-(n+s), \dots, -1, 1, \dots, (n+s)\}$ correspond to $\ell$ and $-i$ to $\ell^0$.
Then $\gn$ contains clauses:
$$
\pp_u(\ab) \rightarrow q^v_i\ \ \ \ \ \ \ \ 
\pp_u(\ab) \rightarrow q^w_{-i}\ \ \ \ \ \ \
\pp_u(\ab) \rightarrow \neg q^u_i\ \ \ \ \ \ \ \ \ 
\pp_u(\ab) \rightarrow \neg q^u_{-i}
$$ 
(these clauses enforce that $\ell$ and $\ell^0$ appear in $D_u, D_v, D_w$ as prescribed by the
resolution rule), 

and for $j \neq i, -i$, $j \in  \{(-(n+s), \dots, -1, 1, \dots, (n+s)\}$
$\gn$ contains further clauses
$$
\pp_u(\ab), q^v_j \rightarrow q^u_j\ \ \ \ \ \ \ \
\pp_u(\ab), q^w_j \rightarrow q^u_j\ \ \ \ \ \ \ \
\pp_u(\ab), q^u_j \rightarrow q^v_j, q^w_j\ \ \ \ \ \ \ \
$$
(these clauses enforce that other literals are passed from $D_v, D_w$ to $D_u$ and
that no other are).

\item $D_u$ was inferred by weakening from $D_v$, $v < u$. Then $\gn$ contains all clauses
$$
\pp_u(\ab), q^v_i \rightarrow q^u_i\ .
$$

\item $D_u$ was inferred as a 1-axiom. Then
$\gn$ contains clause:
$$
\pp_u(\ab) \rightarrow q^u_0 \ .
$$

\end{enumerate}

\item[$\gamma$5.] Finally we add to $\gn$ clauses
$$
\rightarrow \neg q^k_i 
$$
for all $i$, enforcing that $D_k = \emptyset$.

\end{enumerate}
Let us summarize.

\begin{lemma}
For all $n \geq 0$, $s > n$, $k \geq 3s$ the set $\gn$ contains
$O(k^5)$ clauses of width at most $3 + 3 \log k$ and it is not satisfiable.
\end{lemma}

\section{Reductions}

Reflection principles for a proof system Q imply, over an
arbitrary fixed base proof system
satisfying a few technical properties, all Q-provable formulas and only with a polynomial slow-down over
Q. This means that if $\varphi$ has a Q-proof of size $m$ then $\varphi$ can be derived in the base
system from a substitution instance of a reflection principle for Q by a proof of size at most $m^{O(1)}$.
The reader can find all detail in \cite[Sec.9.3]{kniha} but these details are not needed for the
arguments below (although they may help in understanding what is going on).

The set $\gn$ expresses conditions an $\rw$-refutation of some set $\dfc$ would have to satisfy
and hence it is the formula $\neg \bigwedge \gn$ which 
corresponds to reflection principles for ER.
ER-refutation of a set $\Delta$ of clauses amounts to proving formula $\neg \bigwedge \Delta$.
Thus we want derivations (in some base system, here it will be $\rw$) of $\neg \bigwedge \Delta$
from an instance of $\neg \bigwedge \gn$. In the framework of refutation systems this means that
we look for derivations from $\Delta$ of all clauses of a substitution instance of $\gn$. 
In fact, it will be enough to consider $\gz$.

A map $\sigma$ assigning to variables from a set $Y$ constants $0,1$ or disjunctions of
literals corresponding to a set of variables $X$ will be called {\em a clause-substitution}
from $X$ to $Y$,
and the maximal size of a disjunction $\sigma$ assigns is {\em the width} of $\sigma$.

Let $\Gamma, \Delta$ be two sets of clauses in disjoint sets of variables $Y$ and $X$, respectively (to avoid any
confusion when dealing with substitutions). 
We say that {\em $\Delta$ reduces to $\Gamma$ by a clause-substitution $\sigma$} 
iff $\sigma$ is clause-substitution from variables of $\Delta$ to variables of $\Gamma$
such that for each clause $D \in \Gamma$ one of the following cases occurs:

\begin{enumerate}

\item[(a)] $\sigma(D)$ is a 1-axiom,

\item[(b)] $\sigma(D)$ has the form:
\begin{equation} \label{22.8.15b}
\Pi, \bigvee E \rightarrow \bigvee F, \Sigma
\end{equation}
where $E \subseteq F$ are sets of literals. 

\item[(c)] $\sigma(D)$ contains as a subset a clause from $\Delta$,

\end{enumerate} 
Note that in the cases (a) and (b) is $\sigma(D)$ logically valid.

For the construction in the proof of the next theorem it will be handy to use the
following notation. For $z$ a variable and $\ab, \bb \in \at$ put
$$
sel(z,\ab,\bb)
$$
to be the $t$-tuple from $\{0,1,z, \neg z\}^t$ whose $i$-th coordinate is
$$
sel(z,a_i, b_i) := (a_1 \wedge z) \vee (b_i \wedge \neg z) \vee (a_i \wedge b_i)\ .
$$
That is, $sel(z, a_i, b_i)$ is a constant or a literal defined by the following
cases: 
\[ sel(z,a_i, b_i) :=  \left\{ \begin{array}{ll}
 0            &  \mbox{if $a_i = b_i = 0$} \\
 1            &  \mbox{if $a_i= b_i=1$} \\
 z            &  \mbox{if $a_i = 1 \wedge b_i=0$} \\
 \neg z       &  \mbox{if $a_i = 0 \wedge b_i=1$.}
                    \end{array}
            \right. \]

\begin{theorem} \label{main}
Assume $\Delta$ is a set of clauses of width $\le w$ in $n$ variables
that has an ER-refutation $\pi$ with $k(\pi)$ clauses. 

Then for some $k = O(n k((\pi))$ and $s \le k/3$,
$\Delta$ reduces to $\gz$ by a clause-substitution 
of width $\le \max(w,3)$.
\end{theorem}

\prf

Assume $\xx$ are the $n$ variables of $\Delta$. Introducing up to $O(n k(\pi))$ new extensions variables
we may assume the width of $\pi$ is at most $\max(w,3)$. Let $\yy$ be $s$ extensions atoms used in $\pi$.
We may further rearrange the resulting proof
so that the clauses defining the $\yy$ variables are precisely the first $3s$ clauses
and are followed by all $|\Delta|$ clauses from $\Delta$. Let $k = O(n k(\pi))$ be the number
of steps in the resulting ER-refutation and call these steps $D_u$.

Take the set $\gz$ 
and define the following substitution $\sigma$ for its $\pp$- and $\qq$-variables:

\begin{enumerate}

\item For all $\qq$-variables $q^u_i$ with $i \neq 0$ substitute $0$ or $1$, depending on whether
the $\yy$-literal corresponding to $i$ occurs in $D_u$.

\item For all variables $q^u_0$ 
substitute $\bigvee E_u$, where $E_u$ is the set of $\xx$-literals
occurring in $D_u$ together with $1$, if $1 \in D_u$.
(Note that $|E_u| \le w$.)

\item For $\pp$-variables $\pp_u$ with $u= 3r-2, 3-1, 3r$ and $r \le s$ 
define $\sigma$ as follows:
\begin{enumerate}

\item If $D_u$ is one of the three clauses corresponding to an instruction
of the form $y_r := x_j$, put
$$
\sigma(\pp_u)\ :=\ sel(x_j, \ab,\bb)
$$
where $\ab$ and $\bb \in \at$ define the instructions 
$y_r:= 1$ and $y_r:=0$, respectively.

\item Otherwise
substitute for $\pp_u$ the string $\ab^u \in \at$
defining the particular instruction of $\dfc$ in $\pi$.
\end{enumerate}

\item For $\pp_u$ variables with $u = 3s+1, \dots, 3s+|\Delta|$ substitute $\ab \in \at$
defining the clause $D_u$ as being a 1-axiom.

\item For $u = 3s+|\Delta| +1, \dots, k$ consider several cases what to substitute for $\pp_u$:

\begin{enumerate}

\item $D_u$ was inferred as a 1-axiom: substitute for $\pp_u$ as in 
item 4.

\item $D_u$ was derived in $\pi$ by weakening from $D_v$: substitute for $\pp_u$ the $\ab$
specifying this information.

\item $D_u$ was derived by resolution from $D_e, D_f$ resolving variable $y_i$: substitute for $\pp_u$
the $\ab$ specifying this information.

\item As in (c) but the resolved variable was $x_i$. Assume $x_i \in D_e$ and $\neg x_i \in D_f$.
Substitute for $\pp_u$ the expression
$$
sel(x_i, \ab, \bb)
$$
where $\ab, \bb \in \at$ specify that $D_u$ was derived by the weakening from
$D_f$ or $D_e$, respectively.

\end{enumerate}

\end{enumerate}
We need to verify that for every clause $D \in \gz$, $\sigma(D)$ falls under one of the
three cases (a), (b) or (c)
in the definition of reductions by clause-substitutions above. We shall treat
the five groups $\gamma$1 - $\gamma$5 of clauses forming $\gz$ separately.

If $D = \pp_u(\ab) \rightarrow$ belongs to groups $\gamma$1 or $\gamma$3, 
$\sigma(\pp_u(\ab))$ 
contains a false literal and so $\sigma(D)$ is a 1-axiom.

If $D$ is from group $\gamma$2, then $\sigma(D)$ is clearly a 1-axiom by the definition of
$\sigma(q_i^u)$ for all instructions for $y_r$ falling under 3(b) above, i.e. except when it has 
the form $y_r := x_j$. In the latter case the instruction is represented by clauses
$$
\{y_r, \neg x_j\}\ ,\ \ \{\neg y_r, x_j\}\ ,\ \ \{1\}
$$
and $D$ is one of them. The definition of $\sigma$ in 3(a) above
using selection term on $x_j$
yields $\sigma(D)$ which either contains $0$ in the antecedent (and hence $\sigma(D)$ is
a 1-axiom) or one of the literals $x_j, x_j^0$ occurs in both antecedent and succedent
of $\sigma(D)$ and hence it falls under the case (b) of the definition of reductions.  

If $D$ is from group $\gamma$4(a) then 
by item 5(c) of the definition of $\sigma$, $\sigma(D)$ becomes a 1-axiom.
If $D$ is from group $\gamma$4(b) then $\sigma(D)$ falls under the case 
(b) of the definition of reductions:
in particular, for $i=0$, $\sigma(q^v_i)$ is contained in $\sigma(q^u_i)$ (the $E$ and $F$ in
that definition). 
If $D$ is from group $\gamma$4(c) of $\gz$ then $\sigma(D)$ is either a 1-axiom as $\sigma(q^u_0)$ contains constant
$1$ if $D_u$ was a 1-axiom, or it falls under the case (c) of the definition of reductions
as $\sigma(q^u_0)$ is $D_u \in \Delta$ (item 4 of the definition of $\sigma$).

Finally, $D$ from group $\gamma$5 of $\gz$ is trivially turned by $\sigma$ to a 1-axiom.

\qed

We may interpret Theorem \ref{main} as a proof-theoretic reduction: each clause of 
$\sigma(\gz)$ can be derived
from $\Delta$ very easily in any proof system P simulating efficiently the weakening rule and
deriving quickly all 1-axioms and all formulas 
as in (\ref{22.8.15b}) and hence the task to refute $\Delta$ is in P reduced to the task to
refute $\gz$. One can easily list various suitable weak P (e.g. tree-like $R^*(\log)$ or talk about 
$\rw$-derivations of $F,\Pi$ from all $\Sigma, \ell$, $\ell \in E$, in (\ref{22.8.15b})) but it seems
redundant to do so.

Alternatively, we may interpret the theorem
as a reduction between relativized total NP search problems (see e.g. \cite{BCEIP} for definitions). 
That is, given an ER-refutation $\pi$ of $\Delta$ in $n$ variables $\xx$ with
$k(\pi)$ steps,
we have $\gz$ for specific $s, k$ bounded by $O(n k(\pi))$ such that it holds:
\begin{itemize}

\item For any assignment $\alpha$ to variables $\xx$ of $\Delta$ ($\alpha$ is the oracle),
if we know a clause of $\gz$ false under the assignment $\alpha \circ \sigma$ to its variables, we
also know a clause of $\Delta$ false under $\alpha$: $\alpha \circ \sigma(D)$ can only fail 
if it falls under item 2 of the definition
of reductions and hence it contains a clause of $\Delta$ false under $\alpha$.
\end{itemize}
Note that, for a fixed $\pi$ and $\gz$ with parameters determined by it, computing $\sigma$
requires at most $w$ calls to $\alpha$. Hence if $w$ is a constant or at least bounded
by $\log(n |\Delta|)$ the reduction is polynomial time in the sense of \cite{BCEIP}.

It is well-known (see e.g. \cite{kniha,KNT}) that propositional translations of 
a second order $\forall \Sigma^{b}_1(\alpha)$-formula (expressing the totality of a
relativized NP search problem) that is provable in bounded arithmetic
theory $V^1_1$ of Buss \cite{Buss-book}
have polynomial size Extended Frege proofs, i.e. in the refutation 
set-up the corresponding sets of clauses have polynomial size ER-refutations.
This is \cite[Thm.9.1.5]{kniha}, building on earlier results of Cook \cite{Coo75}
and Buss \cite{Buss-book}.

Theorem \ref{main} thus yields the following statement (the definability of the reduction
in $V^0_1$ follows from its explicit nature). 

\begin{corollary} \label{23.8.15a}
Assume that a relativized NP search problem is 
provably total in bounded arithmetic  
theory $V^1_1$.

Then the problem polynomially reduces $\gz$ and the reduction is definable in $V^0_1$. 
\end{corollary}

\section{Total NP search problems $\ig$} \label{22.8.15a}

We shall consider total (non-relativized)
NP search problems given as follows. Let $D(\vv^1, \dots, \vv^t)$ be a circuit
with $t m$ inputs divided into $t$ blocks of size $m$. Such $D$ defines a $t$-ary relation
on $\mm$; as a structure
it may be exponentially large relative to the size of $D$. The general form of search tasks we shall consider
is: Given pair $(1^{(m)}, D)$, find a subset $W \subseteq \mm$ of some specific polynomial
size $m^{O(1)}$ such that the induced substructure is contains a specific configuration
known to exists by a general combinatorial or geometric statement. \cite{Kra-hard} gives
several examples but perhaps the most interesting is when $t = 2$ and
we think of $D$ as defining an undirected graph without loops and 
$W$ either contains a list of $m/2$ vertices from $\mm$ inducing a homogeneous subgraphs
or one or two vertices certifying that $D$ has a lop or is non-symmetric. Ramsey's theorem
$2^m \rightarrow (m/2)^2_2$ guarantees the existence of such a $W$.

We shall use the idea of implicit proofs from \cite{Kra-implicit}, proofs of exponential size
described bit-by-bit 
by a circuit and accompanied by a certificate that the circuit indeed
defines a proof. In particular, a refutation 
of a formula $\phi$ in {\em implicit ER} proof system, denoted iER,
is a pair $(\rho, D)$ such that:
\begin{itemize}

\item $D(u,v)$ is a circuit with two inputs strings $u, v \in \mm$ defining a $2^m \times 2^m$ 0-1 array
which we interpret as describing an ER-refutation of $\phi$ in the same sense as truth assignments to
$\pp$- and $\qq$-variable of $\gz$ talk about a potential ER-refutation,

\item $\rho$ is an ER-proof of the propositional formula formalizing the statement: 

\begin{itemize}

\item {\em $D$ defines a valid ER-refutation of $\phi$}.

\end{itemize}

\end{itemize}
The reader is invited to consult \cite{Kra-implicit} for details of the definition.

The way how we shall use iER was first employed (and justified) in \cite[Thm.5.4]{Kra-hard}.
The idea is simple: we may allow $D$ above to describe not only refutations of polynomial size
formulas (as it was defined in \cite{Kra-implicit})
but of exponential size formulas given themselves by small circuits.

In particular, if $\exists y (|y| \le |x|^c) \varphi(x,y)$ is a $\Sigma^b_1$-formula
with $\varphi \in \Sigma^b_0$, then the sentence 
$\forall x (|x|=n) \exists y (|y| \le |x|^c) \varphi(x,y)$
is true iff the set
\begin{equation}\label{29.8.15a}
\neg \varphi(x, w)\ , \mbox{ all $w$ such that }\ |w| \le n^{c} 
\end{equation}
is not satisfiable by any $x \in \nn$, and hence it is refutable (in ER, in particular).
Set (\ref{29.8.15a}) has exponential size but it can be easily generated by a size $n^{O(c)}$
circuit from $w$'s. We use these ideas as follows.

\bigskip

The {\bf $\igm$ NP search problem}, the instance of $\ig$ for parameter $m$, is defined as follows:

\begin{enumerate}

\item The input is pair $(1^{(m)}, D)$ with $D(x,y)$ a size $m^2$ circuit with 
$2 \cdot m$ inputs.

\item Interpret $D$ as defining an evaluation to $\pp$ and $\qq$-variables of $\gz$
where $k = 2^m$ and $s = k^{1/2}$.

[There are $2s+1$ variables $q^u_i$ and $3 \log k = 3m$ variables in $\pp_u$, all $u \le k = 2^m$,
so $D$ has enough input bits to define a $0-1$ array of bits evaluating all these variables.]

\item Output: find a clause of $\gz$ false under the evaluation. 

[There are $O(k^5) = O(2^{5m})$ of possible outcomes so the output is $\le 5m $ bits.]
\end{enumerate}
The parameters are fixed at $|D|= m^2$ and $s = 2^{m/2}$ in order to reduce the number of
parameters in the problem. Modifying $m$ linearly allows to accommodate arbitrary polynomial 
relations among $k, s$ and $\log k, |D|$.

\medskip

We state and prove the next theorem using the ideas and referring to
facts about the concepts described above; all
details for these facts can be found in \cite{Kra-implicit,Kra-hard}
at the specifically cited places.

\begin{theorem} \label{23.8.15b}
Assume an NP search problem is provably total in theory $V^1_2$. 
Then the problem can be polynomially
reduced to $\ig$. The reduction is definable in $S^1_2$.

Moreover, $\ig$ is itself provably total in $V^1_2$.
\end{theorem}

\prf

Let $\exists y (|y| \le |x|^c) \varphi(x,y)$ be a $\Sigma^b_1$-formula 
with $\varphi \in \Sigma^b_0$ such that $V^1_2$ proves
$$
\forall x \exists y (|y| \le |x|^c) \varphi(x,y)\ .
$$ 
In particular, $|y| \le |x|^c \wedge \varphi(x,y)$
defines a total NP search problem.

By the construction underlying \cite[Thm.2.1]{Kra-implicit}, as shown in the proof of \cite[Thm.5.4]{Kra-hard},
there exists an iER refutation $(\rho, B)$ of formulas from (\ref{29.8.15a}) above
expressed as a set of $2^{O(n^c)}$ clauses of width $\le w = n^{O(c)}$ such that:

\begin{itemize}

\item $B(i,j)$ is 
a size $n^{O(c)}$ circuit $B(i,j)$ with $2 \cdot n^{O(c)}$ inputs 
describing an ER-refutation $\pi$ of (\ref{29.8.15a}),

\item circuit $B$ is definable in $S^1_2$ from $1^{(n)}$ and $S^1_2$ proves that $B$ 
defines an ER-refutation of the set (\ref{29.8.15a}) ($\rho$ plays a role in this).

\end{itemize}
Use $\pi$ for the definition of a clause-substitution $\sigma$ as in the proof of Theorem \ref{main}
but whenever we need a bit of $\pi$ we compute it by circuit $B$. The substitution has width $\le w
= n^{O(c)}$ and so the reduction so obtained is a polynomial reduction of the search problem
$\exists y (|y| \le |x|^c) \varphi(x,y)$ for $x$ of length $|x|=n$.

\smallskip

The second statement follows as $V^1_2$ proves the soundness iEF (= iER) proofs,
cf. \cite[Thm.2.1]{Kra-implicit}.

\qed

One can generalize this construction to stronger theories as follows. 
In \cite{Kra-implicit}
we used the characterization from \cite{Kra-exp}
of bounded first-order consequences of $V^1_2$ as those of formal system 
$S^1_2 + \mbox{1-Exp}$ : $\delta(x)$ is provable in this system iff $S^1_2$ proves
$$
t(x) \le |y| \rightarrow \delta(x)
$$
for some term $t(x)$. The intuition is that while ER corresponds to $S^1_2$,
iER corresponds to adding 1-Exp and that corresponds to extending first-order
$S^1_2$ to second-order $V^1_2$. The construction in \cite{Kra-implicit,Kra-hard}
works also for $S^1_2 + \mbox{2-Exp}$ (and third-order extension of $S^1_2$) and i(iER),
and higher iterates, as pointed out in \cite[Sec.4]{Kra-implicit}. 
In general, if a theory T corresponds to a proof system P then 
iP corresponds to T + 1-Exp and one may try to define NP
search problems analogous to $\ig$ where $B$ is assumed to describe a P-refutation.
It is a challenge to describe this construction in a direct, combinatorially
transparent, way.

\section{Concluding remarks}

One can restrict circuits $C$ that can be used in $\dfc$ to a class of circuits
and it is clear from the construction that
taking for these classes $NC^1, AC^0$ or $AC^0(2)$ yield statements analogous to Theorem 
\ref{main}
for Frege system and constant depth Frege system in the DeMorgan language,
and constant depth Frege systems
in DeMorgan language augmented by the parity connective, respectively.
Similarly, Corollary \ref{23.8.15a} and Theorem \ref{23.8.15b}
can be analogously derived for theories corresponding to those proof systems, cf.
\cite{kniha,CN}.

For a given $s \geq 1$ and $k \geq 3s$ we can define the following
random process yielding a set of clauses ($\mathbf r$ are the random bits used):

\begin{enumerate}

\item Pick $s$ instructions for computing variables $\yy$ defining a circuit
$C_{\mathbf r}$ without variables: the 
instruction for $y_i$ is picked uniformly at random from all legal instructions
for $y_i$,

\item substitute in $\gz$ for all variables $\pp_u$ and $q^u_i$ with $u\le 3s$
the bits defining clauses of $\dfz$ corresponding to $C_{\mathbf r}$ chosen in step 1.

\end{enumerate}
Let us denote the random set of clauses so constructed by $\gz(C_{\mathbf r})$; it is always
unsatisfiable. The following seem to be interesting open problems:
\begin{enumerate}

\item Is it true that with a high probability over $\mathbf r$ the set $\gz(C_{\mathbf r})$
requires long refutations in any proof system not simulating ER?

\item Is it true that $\gz$ can be reduced by a clause-substitution 
to a problem 
$\Gamma(0,s',k')(C)$ for some $k'$ polynomially
bounded in $k$ and some {\em specific} size $s'$ circuit $C$?

\end{enumerate}
If the first question had an affirmative answer sets
$\gz(C_{\mathbf r})$ would provide an easy to compute source of hard
formulas that are always unsatisfiable (other proposed constructions yield sets
unsatisfiable with a high probability but not always).

\bigskip

\noindent
{\large {\bf Acknowledgements:}}

\medskip

The constructions described in this paper were developed as a part of an
investigation into model-theoretic constructions described in \cite{Kra-forcing}.
It was after discussions with S.~Buss in July 2015 about his current work \cite{Bec-Bus}
with A.~Beckmann that I realized that the constructions could be of an independent
interest. I also thank S.~Buss for comments on the draft of this paper and to 
N.~Thapen for discussions about the topic.

\bigskip
\noindent
{\bf Mailing address:}

Department of Algebra

Faculty of Mathematics and Physics

Charles University

Sokolovsk\' a 83, Prague 8, CZ - 186 75

The Czech Republic

{\tt krajicek@karlin.mff.cuni.cz}

\end{document}